\documentclass[11pt,a4paper]{amsart}
\usepackage{epsfig}

\def\N{\mathbb{N}}
\def\Z{\mathbb{Z}}
\def\R{\mathbb{R}}
\def\Q{\mathbb{Q}}
\def\C{\mathbb{C}}
\def\H{\mathbb{H}}

\def\Re{{\rm Re}}
\def\Im{{\rm Im}}
\def\eps{\varepsilon}
\def\Cc{\widehat{{\C}}}
\def\ol{\overline}
\def\tilde{\widetilde}
\def\sm{\setminus}
\def\id{{\rm id}}
\def\a{{\alpha}}
\def\e{{\varepsilon}}

\def\cal{\mathcal}

\newtheorem{theorem}{Theorem}[section]
\newtheorem{proposition}[theorem]{Proposition}
\newtheorem{lemma}[theorem]{Lemma}
\newtheorem{definition}[theorem]{Definition}
\newtheorem{corollary}[theorem]{Corollary}
\def\proofof #1 {\par\medskip\noindent {\it Proof of #1. }}
\def\remark{\par\medskip \noindent {\sc Remark. }}

\title{Virtual Immediate Basins of Newton Maps and Asymptotic Values}
\author{Xavier Buff}
\address{Universit\'e Paul Sabatier, 118 route de Narbonne, 31062 Toulouse, France}
\email{buff@picard.ups-tlse.fr}
\author{Johannes R\"uckert}
\address{International University Bremen, Campus Ring 12, 28759 Bremen, Germany}
\email{j.rueckert@iu-bremen.de} \subjclass[2000]{30D05, 30D99,
37F10, 37F75, 49M15} \keywords{Newton map, virtual immediate basin,
asymptotic value, Baker domain}

\begin{document}
\renewcommand\labelenumi{\arabic{enumi}.}

\begin{abstract}
Newton's root finding method applied to a (transcendental) entire
function $f:\C\to\C$ is the iteration of a meromorphic function
$N_f$. It is well known that if for some starting value $z_0$,
Newton's method converges to a point $\xi\in\C$, then $f$ has a root
at $\xi$. We show that in many cases, if an orbit converges to
$\xi=\infty$ for Newton's method, then $f$ has a `virtual root' at
$\infty$. More precisely, we show that if $N_f$ has an invariant
Baker domain that satisfies some mild assumptions, then $0$ is an
asymptotic value for $f$.

Conversely, we show that if $f$ has an asymptotic value of
logarithmic type at $0$, then the singularity over $0$ is contained
in an invariant Baker domain of $N_f$, which we call a {\em virtual
immediate basin}. We show by way of counterexamples that this is not
true for more general types of singularities.
\end{abstract}

\maketitle

\section{Introduction}
Let $f:\C\to\C$ be an entire function. Newton's root finding method
for $f$ is implemented by iterating the associated {\em Newton map}
\[
    N_f:\C\to\Cc, \;z\mapsto z-\frac{f(z)}{f'(z)}\;.
\]
It is well known that $\xi\in\C$ is a fixed point of $N_f$ if and
only if $f(\xi)=0$. Furthermore, every finite fixed point $\xi$ of
$N_f$ is attracting, so it has an invariant neighborhood on which
$N_f$-orbits converge locally uniformly to $\xi$. In 2003, Douady
raised the following question: if there exists a {\em virtual
immediate basin} (an invariant, unbounded domain on which
$N_f$-orbits converge locally uniformly to $\infty$), does this
imply that $\infty$ is a `virtual root' of $f$, in other words, does
this imply that $0$ is an {\em asymptotic value} of $f$? In this
paper, we give a condition under which this is true. A recent result
of Bergweiler, Drasin and Langley \cite{BDL} implies that the condition is
sharp when the Julia set of Newton maps is connected. Conversely, we
show that if $f$ has a singularity of logarithmic type over $0$,
then this singularity is contained in a virtual immediate basin of
$N_f$; if it is not of logarithmic type, then we provide
counterexamples.

The dynamics of $N_f$
partitions the Riemann sphere $\Cc$ into two completely invariant
parts: the open {\em Fatou set} of all points at which the iterates
$\{N_f^{\circ n}\}_{n=0}^{\infty}$ are defined and form a normal
family in the sense of Montel, and its complementary {\em Julia set}
that contains the backward orbit of $\infty$; see \cite{Bergweiler,
Milnor} for an introduction to these concepts. Note that starting
values in the Julia set will never converge to an attracting fixed
point of $N_f$.

A component of the Fatou set of $N_f$ for which no point converges
to a root of $f$ under iteration is either {\em wandering} or will
eventually land on a cycle of {\em B\"ottcher domains}, {\em Leau
domains}, {\em Siegel disks}, {\em Herman rings} or {\em Baker
domains} (compare \cite[Theorem 6]{Bergweiler}).

The possibilities become much more restricted when considering an
{\em invariant} component $U$ of the Fatou set, so that
$N_f(U)\subset U$. In this case, it follows from Proposition
\ref{Prop_NewtonMaps} that $U$ either contains a root of $f$, or is
an invariant Herman ring or Baker domain.

Shishikura \cite{Shishikura} has shown that if $N_f$ is rational,
then its Julia set is connected (see Proposition
\ref{Prop_RationalNewton} for a characterization of rational Newton
maps). It is conjectured that Shishikura's result can be extended to
all Newton maps of entire functions. If this is true, an invariant
Fatou component of $N_f$ either contains a root of $f$ or is a
virtual immediate basin (see Section \ref{Sec_Baker} for the precise
definition).

In this paper, we continue the analysis of virtual immediate basins
in \cite{MS} and \cite{RS}. We prove that if $f$ has a logarithmic
singularity over $0$, then $N_f$ has a virtual immediate basin (in
1994, Bergweiler, von Haeseler, Kriete, Meier and Terglane
investigated a class of functions $f$ that tend to $0$ in a sector
and showed that a right end of this sector is contained in a Baker
domain of $N_f$ \cite[Theorem 3.3]{BHK}).

For non-logarithmic singularities over $0$, we give examples of
functions whose Newton maps do not have a virtual immediate basin
associated to these singularities.

Furthermore, we show that there are three classes of virtual
immediate basins for $N_f$, two of which induce an asymptotic value
at $0$ for $f$. For the third class, this statement requires an
additional assumption, without which it is false. Every such virtual
immediate basin even has an open subset of starting values $z_0$
such that as $z_n=N_f^{\circ n}(z_0)\to\infty$, $f(z_n)\to 0$.

\medskip
Our paper is structured as follows: in Section \ref{Sec_Baker}, we
give a precise definition of virtual immediate basins and state
several of their properties. In Section \ref{Sec_Singularities}, we
recall some fundamental notions concerning singular values. In
Section \ref{Sec_Main}, we prove that a logarithmic singularity over
$0$ for $f$ induces a virtual immediate basin for $N_f$, while the
counterexamples for direct singularities are treated in Section
\ref{Sec_Example}. The converse theorem is stated and proved in
Section \ref{Sec_Converse}. The underlying idea of the proof is to
compare iterates of the Newton map $\displaystyle
N_f=\id-\frac{f}{f'}$ to the time 1 flow of $\displaystyle \dot
z=-\frac{f(z)}{f'(z)}$.

\section{Virtual Immediate Basins}
\label{Sec_Baker} The concept of a {\em virtual immediate basin} was
introduced in \cite{MS} to explain the behavior of Newton maps
between different accesses to $\infty$ of an immediate basin.
Examples of Newton maps having virtual immediate basins can be found
in \cite{MS, RS}; these example are discussed in detail in
\cite{Mayer}. The name was chosen to suggest that these domains
behave in many ways similar to immediate basins.

The following proposition characterizes Newton maps of entire
functions.
\begin{proposition}[Newton Maps]
\label{Prop_NewtonMaps} {\em \cite[Proposition 2.8]{RS}.} Let
$N:\C\to\Cc$ be a meromorphic function. It is the Newton map of an
entire function $f:\C\to\C$ if and only if for each fixed point
$N(\xi)=\xi\in\C$, there exists a natural number $m>0$ such that
$N'(\xi)=\frac{m-1}{m}<1$. In this case, there exists $c\neq 0$ such
that
\[
    f = c\cdot \exp\left(\int \frac{d\zeta}{\zeta-N(\zeta)}\right)\;\;.
\]
\qed
\end{proposition}

Note that while all definitions in this section are written in terms
of Newton maps, they make sense for arbitrary meromorphic functions.
\begin{definition}[Immediate Basin]
\label{Def_ImmediateBasin} Let $N_f$ be a Newton map. If $\xi$ is an
attracting fixed point of $N_f$, we call the open set
\[
    \{z\in\C\,:\,\lim_{n\to\infty}N_f^{\circ n}(z)=\xi\}
\]
its {\em basin (of attraction)}. The component of the basin that
contains $\xi$ is called its {\em immediate basin} and denoted
$U_{\xi}$.
\end{definition}
For the definition of virtual immediate basins, we need the
following concept.
\begin{definition}[Absorbing Set]
\label{Def_AbsorbingSet} Let $V$ be an $N_f$-invariant domain. A
connected and simply connected open set $A\subset V$ is called a {\em weakly absorbing set for
$V$} if $N_f(A)\subset A$ and for each compact
$K\subset V$, there exists $k\in\N$ such that $N_f^{\circ
k}(K)\subset A$. 

We call $A$ an {\em absorbing set} if it is weakly absorbing and additionally satisfies $N_f(\ol{A})\subset A$, where the closure is taken in $\C$.
\end{definition}
\begin{definition}[Virtual Immediate Basin]
\label{Def_VirtualBasin} A domain $V\subset\C$ is called a {\em
virtual immediate basin} for $N_f$ if it is maximal (among domains in $\C$) with respect to
the following conditions:
\begin{enumerate}
\item for every $z\in V$, $\lim_{n\to\infty} N_f^{\circ n}(z)=\infty$;
\item $V$ contains an absorbing set.
\end{enumerate}
\end{definition}
Every virtual immediate basin is unbounded, invariant and simply
connected \cite[Theorem 3.4]{MS}. Since Newton maps of polynomials
have a repelling fixed point at $\infty$, virtual immediate basins
can appear only for Newton maps of transcendental functions.
\begin{proposition}[Rational Newton Map]
\label{Prop_RationalNewton} {\em \cite[Proposition 2.11]{RS}.} Let
$f:\C\to\C$ be an entire function. Its Newton map $N_f$ is rational
if and only if there exist polynomials $p,q$ such that $f=p\cdot
e^q$. In this case, $\infty$ is a repelling or parabolic fixed
point.

More precisely, let $m:= \deg p$ and $n:=\deg q$. If $n=0$ and
$m\geq 2$, then $\infty$ is repelling with multiplier
$\frac{m}{m-1}$. If $n>0$, then $\infty$ is parabolic with
multiplier $+1$ and multiplicity $n+1\geq 2$. \qed
\end{proposition}
In the following, let $f$ be a transcendental entire function. If
$N_f$ is rational, then it has virtual immediate basins which are
the attracting petals of the parabolic fixed point at $\infty$ (see
\cite[Theorem 10.5]{Milnor}). If $N_f$ is transcendental
meromorphic, then any virtual immediate basin is (contained in) an
invariant Baker domain.
\begin{definition}[Baker Domain]
\label{Def_BakerDomain} Let $B$ be an invariant component of the
Fatou set of $N_f$. If $\lim_{n\to\infty}N_f^{\circ
n}(z)=\infty\in\partial B$ for all $z\in B$ and $N_f$ has an
essential singularity at $\infty$, then we call $B$ a {\em Baker
domain} of $N_f$.
\end{definition}
If $B$ is a simply connected Baker domain, it contains a weakly absorbing
set $A$ by a result of Cowen \cite[Theorem 3.2]{Cowen}. Using Cowen's work, it is easy to find an absorbing subset of $A$, hence $B$ is a
virtual immediate basin. Moreover, Cowen's result implies that there
are three dynamically defined classes of virtual immediate basins.
The following notations are based on \cite{Koenig} and
\cite{FagellaBaranski}.
\begin{definition}[Conformal Conjugacy]
\label{Def_ConformalConjugacy} Let $V$ be a virtual immediate basin
of $N_f$ and define $T(z)=z+1$. If there exists a weakly absorbing set $A$
for $V$, a $T$-invariant domain $\Omega\subset\C$ and a holomorphic
map $\phi:V\to\Omega$ such that
\begin{equation*}
    \phi\circ N_f(z) = T\circ \phi(z)
\end{equation*}
for all $z\in V$, $\phi$ is univalent on $A$ and
$\phi(A)\subset\Omega$ is a weakly absorbing set for $T|_{\Omega}$, then we call the
triple $(\Omega,\phi,T)$ a {\em conformal conjugacy} for $N_f$ on $V$.
\end{definition}
\begin{definition}[Types of Virtual Immediate Basins]
\label{Def_Types} Let $V$ be a virtual immediate basin of $N_f$. We
say that $V$ is {\em parabolic of type I} if it has a conformal
conjugacy $(\Omega,\phi,T)$ such that $\Omega=\C$. It
is {\em parabolic of type II} if there exists a conjugacy such that
$\Omega$ is an upper or lower half-plane and {\em hyperbolic} with
constant $h$ if there exists $h>0$ such that $\Omega$ is the strip
\[
    S_h:=\{z\in\C\,:\,|\Im(z)|<h\}\;.
\]
\end{definition}
\begin{theorem}[Classification of Virtual Immediate Basins]
\label{Thm_Classification} {\em \cite[Theorem 3.2]{Cowen}.} Every
virtual immediate basin $V$ has a conformal conjugacy and is of
exactly one of the three types defined above. If $V$ is hyperbolic,
the constant $h$ is uniquely defined. \qed
\end{theorem}
\remark We believe that any Baker domain of a Newton map is simply
connected; if this were proved, the notion of a virtual immediate
basin would simply stand for either an attracting petal or a Baker
domain, depending on whether the map under consideration is rational
or not.

\section{Asymptotic Values}
\label{Sec_Singularities} We recall several important definitions
concerning the singular values of a meromorphic map. Singular values
play an important role in iteration theory, because their orbits
determine the dynamics of a map in many ways.

We denote by $B_r(z)$ the open disk of radius $r>0$ around $z\in\C$.
In this section, let $g:\C\to\Cc$ be a meromorphic function.

\begin{definition}[Regular and Singular Value]
Let $a\in\C$ and assume that for $r>0$, $U_r$ is a connected
component of $g^{-1}(B_r(a))$ such that $U_{r_1}\subset U_{r_2}$ if
$r_1<r_2$.\footnote{The function $U:r\mapsto U_r$ is completely
determined by its germ at $0$. Since $\bigcap_{r>0} U_r$ is
connected, the intersection contains at most one point.} We have the
following two cases:
\begin{enumerate}
\item
If $\bigcap_{r>0}U_r=\{z\}$ for some $z\in\C$, then $g(z)=a$. If
$g'(z)\neq 0$, then we call $z$ a {\em regular point} of $g$. If
$g'(z)=0$, then $z$ is called a {\em critical point} and $a$ a {\em
critical value}. In this case, we say that the critical point $z$
{\em lies over} $a$.
\item
If $\bigcap_{r>0}U_r=\emptyset$, then we say that $U:r\mapsto U_r$
defines a {\em singularity of $f^{-1}$} and we call $a$ an {\em
asymptotic value}. For simplicity, we call $U$ a {\em singularity}
and say it {\em lies over} $a$.
\end{enumerate}
A {\em singular value} is an asymptotic or critical value. If no
singularities or critical points lie over a point, we call it a {\em
regular value}.
\end{definition}
Note that there can be many different singularities as well as
regular or critical points over any given point $a\in\C$.

For a rational map, all singular values are critical values.
Asymptotic values of transcendental maps have a well-known
characterization via paths.
\begin{lemma}[Asymptotic Path]
A point $a\in\Cc$ is an asymptotic value of $g$ if and only if there
exists a path $\Gamma:(0,\infty)\to\C$ with
$\lim_{t\to\infty}\Gamma(t)=\infty$ such that
$\lim_{t\to\infty}g(\Gamma(t))=a$. \qed
\end{lemma}
We call $\Gamma$ an {\em asymptotic path} of $a$. We follow
\cite{BergweilerEremenko} in the classification of asymptotic
values.
\begin{definition}[Direct, Indirect and Logarithmic Singularity]
\label{Def_LogSing} Let $U$ be a singularity of $g^{-1}$ lying over
$a\in\C$.

If $a\not\in g(U_r)$ for some $r>0$, then we call $U$ a {\em direct}
singularity. Otherwise, $U$ is called an {\em indirect} singularity.

A direct singularity $U$ over $a$ is called {\em logarithmic} if $g:
U_r\to B_r(a)\setminus\{a\}$ is a universal covering map for all
sufficiently small $r$.
\end{definition}
 As an example, the positive real axis is an asymptotic path of $0$ for the
 map $z\mapsto \sin(z)/z$. Since its image assumes this value infinitely many
 times, it is contained in an indirect singularity over $0$.
For $z\mapsto \exp z$, any left half plane is a logarithmic
singularity over $0$.

\section{A Criterion for Virtual Immediate Basins}
\label{Sec_Main}
Our first result is the following.
\begin{theorem}[Logarithmic Singularity Implies Virtual Immediate Basin]
\label{Thm_Logarithmic} Let $f:\C\to\C$ be an entire function with a
logarithmic singularity $U$ over $0$. Then there exists $r_0>0$ such
that $U_{r_0}$ is an absorbing set for a parabolic virtual immediate
basin of type I for $N_f$.
\end{theorem}
Note that if $U$ is an indirect singularity, each $U_r$ contains
infinitely many roots of $f$ and hence infinitely many attracting
fixed points of $N_f$. Therefore, $U_r$ cannot be part of a virtual
immediate basin. In Section \ref{Sec_Example}, we show that there
exist functions $f:\C\to\C$ with a direct singularity $U$ over $0$
which does not induce a virtual immediate basin for $N_f$.

\begin{proof}
The idea is to compare the iterates of $N_f$ to the time 1 flow of
the differential equation $\dot z = -\frac{f(z)}{f'(z)}$. If $r$ is
small enough, this flow sends $U_r$ isomorphically to $U_{r/e}$. We
will see that for $r$ small enough, $N_f$ maps $U_r$ univalently
into itself and is an absorbing set for a virtual immediate basin of
$N_f$.

First, let $r>0$ be small enough so that $f:U_r\to B_r(0)\sm\{0\}$
is a universal covering. Set $\eta:=-\log r$ and
$\H_\eta:=\{w\in\C\,:\, \Re(w)> \eta\}$. Since $e^{ -\id}:  \H_\eta
\to B_r(0)\sm\{0\}$ is also a universal covering, the map
$-\log(f):U_r\to\H_\eta$ is biholomorphic with inverse $\psi:\H_\eta
\to U_r$ (see Figure \ref{Fig_1}). With this, we get $\log (f(
\psi(w)))=-w$ for $w\in\H_\eta$.
\begin{figure}[hbt]
\centerline{
\begin{picture}(0,0)%
\includegraphics{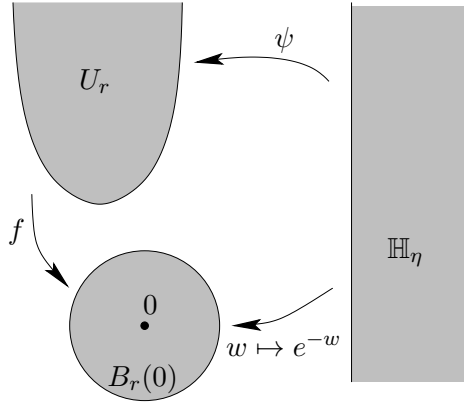}%
\end{picture}%
\setlength{\unitlength}{2072sp}
\begin{picture}(5356,4787)(3196,-4341)
\put(7741,-2626){$\H_\eta$} \put(4051,-601){$U_r$}
\put(6391,-106){$\psi$} \put(5811,-3796){$w\mapsto e^{-w}$}
\put(4391,-4156){$B_r(0)$} \put(4816,-3301){$0$}
\put(3196,-2401){$f$}
\end{picture}}
\caption{\label{Fig_1} If $f:U_r\to B_r(0)\sm\{0\}$ is a universal
covering, there exists a biholomorphic map $\psi: \H_\eta\to U_r$.}
\end{figure}

Taking derivatives yields
\[
    \frac{f'(\psi(w))}{f(\psi(w))} \cdot \psi'(w) = -1; \quad\text{hence}\quad
\psi'(w) = -\frac{f(\psi(w))}{f'(\psi(w))}\;.
\]
In other words, $\psi$ is a solution of $\dot z=
-\frac{f(z)}{f'(z)}$ and following the flow during time 1 maps
$U_r=\psi(\H_\eta)$ to $U_{r/e}=\psi(\H_{\eta+1})$.

We now want to compare $N_f$ to the time 1 flow of $\dot z=
-\frac{f(z)}{f'(z)}$. We will do the comparison in time space: we
will show that if $z=\psi(w)$ with $\Re(w)$ large enough, then
$N_f(z) = \psi(w')$ with $w'$ close to $w+1$. More precisely, we
have the following lemma.

\begin{lemma}
\label{Lem_GconjNf} There exists $\eta_0>\eta$ and a holomorphic map
$G:\H_{\eta_0}\to \H_{\eta_0+1/2}$ such that for all $w\in
\H_{\eta_0}$, we have
\[N_f\circ \psi(w) = \psi\circ G(w)\quad\text{and}\quad
|G(w)-(w+1)|<\frac{1}{2}.\]
\end{lemma}

The proof of Theorem \ref{Thm_Logarithmic} is then easily completed.
Indeed, set $V_0:=\psi(\H_{\eta_0})=U_{r_0}$ with $r_0=e^{-\eta_0}$
and let $V_{n+1}$ be the component of $N_f^{-1}(V_{n})$ that
contains $V_0$. Since all points in $\H_{\eta_0}$ converge to
$\infty$ under iteration of $G$ (the real part increases by at least
$1/2$ in each step), we conclude that $V:=\bigcup_{n\in\N} V_n$ is a
virtual immediate basin of $N_f$ with absorbing set $V_0$.

Let us now prove Lemma \ref{Lem_GconjNf}. Note that
\[
    N_f(\psi(w))=\psi(w) - \frac{f(\psi(w))}{f'(\psi(w))}=\psi(w) + \psi'(w)\;.
\]
Thus, it is equivalent to prove that there exists $\eta_0>\eta$ and
a holomorphic map $G:\H_{\eta_0}\to \H_{\eta_0+1/2}$ such that for
all $w\in \H_{\eta_0}$, we have
\begin{equation}
\label{Eqn_ToShow} \psi(w) + \psi'(w)  = \psi(G(w))
\quad\text{and}\quad |G(w) - (w+1)| < \frac{1}{2}\;.
\end{equation}

Given $w\in \H_{\eta+2}$, define functions $g,h: B_2(w)\to \C$ by
\[
g:\zeta \mapsto \frac{\psi(\zeta)-\psi(w)-\psi'(w)}{\psi'(w)}
\quad\text{and} \quad h:\zeta\mapsto  \zeta-(w+1)\;.
\]
Since $g$ and $h$ satisfy $g(w)=h(w)=-1$, $g'(w)=h'(w)=1$ and can
both be extended to all of $\H_\eta$ as univalent maps, by Koebe's
distortion theorem there exists $\eta_0>\eta+2$ such that for every
$w\in \H_{\eta_0}$ and every $\zeta\in B_2(w)$,
$|g(\zeta)-h(\zeta)|<1/4$.

Clearly, $h(w+1)=0$. Note that $|h(\zeta)|=1/2>|g(\zeta)-h(\zeta)|$
when $\zeta$ belongs to the circle $\partial B_{1/2}(w+1)$. By
Rouch{\'e}'s theorem, the map $g$ has a (unique) root $\xi_w \in
B_{1/2}(w+1)$. It is now easy to see that the map $G:\H_{\eta_0}\to
\C$ defined by $G(w)=\xi_w$ satisfies equations (\ref{Eqn_ToShow}).
\end{proof}

\section{A Direct Singularity Counterexample}
\label{Sec_Example}

In this section, we will exhibit examples of entire functions with
direct singularities over $0$ that do not induce Baker domains of
the associated Newton maps. This shows that Theorem
\ref{Thm_Logarithmic} cannot be improved much further; $0$ is an
omitted value in all examples, so that a generalization is not even
possible to this case. We will only treat the first example in full
detail.

For $\a\in \left]0,+\infty\right[$, consider the entire function
$f_\a$ defined by
\[
f_\a(Z) = \exp\left(-\frac{1}{\a}\left(Z+\frac{1}{2i\pi} e^{2i\pi
Z}\right)\right) \;.
\]
The function $f_\a$ has infinitely many singularities over $0$ which
are necessarily direct since $f_\a$ does not vanish. We have two
kinds of asymptotic paths:
\begin{enumerate}
\item \label{sing_type1}for $k\in \Z$, as $t\in \R\to +\infty$, $f_\a(k+\frac{1}{4}-it)\to 0$;
\item \label{sing_type2}as $t\in \R\to +\infty$, $f_\a(t)\to 0$.
\end{enumerate}
The singularities of the first kind are of logarithmic type. Thus,
each one induces a Baker domain of parabolic type I for the Newton
map
\[
        N_\a(Z) = Z+\frac{\a}{1+e^{2i\pi Z}} \;.
\]
The singularity of the second kind is not of logarithmic type and
contains infinitely many critical points of $f$. We will see that
for some values of $\a$, it does not induce a Baker domain for
$N_\a$.

More precisely, observe that $N_\a(Z+1)=N_\a(Z)+1$. It follows that
we can study the dynamics of $N_\a$ modulo $1$. In other words, we
have
\[
        e^{2i\pi N_\a(Z)} = g_\a\left(e^{2i\pi Z}\right)\quad\text{with}\quad
g_\a(z) = ze^{2i\pi\a/(1+z)} \;.
\]
The map $g_\a$ has a fixed point with multiplier $e^{2i\pi \a}$ at
$z=0$, a fixed point with multiplier $1$ at $z=\infty$ and an
essential singularity at $z=-1$.

Let ${\cal F}(N_\a)$ and ${\cal F}(g_\a)$ be the Fatou sets of
$N_\a$ and $g_\a$ and let $\pi: \C\to \C^*$ be the universal
covering $\pi:Z\mapsto z=e^{2i\pi Z}.$ We claim that
\[
{\cal F}(N_\a) = \pi^{-1}\bigl({\cal F}(g_\a)\bigr) \; .
\]
It is easy to see that $\pi^{-1}\bigl({\cal F}(g_\a)\bigr)\subset
{\cal F}(N_\a)$ (see for example \cite{Bergweiler2}). The inclusion
${\cal F}(N_\a)\subset \pi^{-1}\bigl({\cal F}(g_\a)\bigr)$ is less
immediate. One may argue as follows. Assume $z_0=\pi(Z_0)\notin
{\cal F}(g_\a)$. Then, $z_0$ lies in the closure of the set of
iterated $g_\a$-preimages of $-1$ (otherwise, the family of iterates
of $g_\a$ would be well defined near $z_0$ and avoid the infinite
set $g_\a^{-1}(\{-1\})$, thus it would be normal). It follows that
any neighborhood of $Z_0$ contains a preimage of a pole of $N_\a$.
Thus, $Z_0\notin {\cal F}(N_\a)$.

As $z\to \infty$, we have
\[
        g_\a(z) = z + 1+ \frac{2i\pi\a}{z}+o(1/z)\;.
\]
Thus, the parabolic fixed point at $\infty$ has multiplicity $2$. It
has a single attracting direction along the positive real axis. The
full preimage of its parabolic basin under the map $e^{2i\pi Z}$ is
the union of the Baker domains of $N_\a$ induced by the
singularities of $f_\a$ of the first kind. The map $g_\a$ has
exactly two critical points: the solutions to $(1+z)^2-2 i \pi\a
z=0$.

Conjugating with $z\mapsto w=1/(z+1)$, we may put the singularity at
$\infty$ and the fixed points at $0$ and $1$. The map $g_\a$ is thus
conjugate to the meromorphic function
\[
        h_\a(w) = \frac{w}{w+(1-w)e^{2i\pi \a w}} \;.
\]
The map $h_\a$ has growth order $1$ and two critical points. Thus,
it has at most $2$ asymptotic values by \cite[Corollary
3]{BergweilerEremenko}. But as $t\in \R\to +\infty$, $h_\a(it)\to
0$ and $h_\a(-it)\to 1$. Thus, $h_\a$ has exactly $2$ (fixed)
asymptotic values and $2$ critical values and is therefore a finite
type map. It is well known that finite type meromorphic functions
have neither wandering domains nor Baker domains \cite{BakerKotusLu,
RipponStallard}.

The map $h_{\a}$ has a fully invariant parabolic point at $0$ and
for suitably chosen $\a$, the fixed point at $1$ is Cremer (in
analogy to \cite[Theorem 11.13]{Milnor}). We want to prove that in
this case, the Fatou set of $h_{\a}$ consists of the parabolic basin
at $0$ and its preimage components. We deduce that then, the Fatou
set of $g_\a$ is equal to the parabolic basin of $\infty$ and its
preimage components. Thus, every Fatou component of $N_\a$ maps
after finitely many iterations into one of the invariant Baker
domains induced by the first kind of singularities of $f_\a$. There
is no Fatou component associated to the second kind of singularity
of $f_\a$.

So it remains to show that $h_{\a}$ has no additional non-repelling
periodic points nor Herman rings. While both claims follow directly
from Epstein's version of the Fatou-Shishikura inequality for finite
type maps \cite{Epstein1, Epstein2, Epstein3}, we provide a version
of Epstein's proof that is sufficient for our purposes; we treat
Herman rings separately in Lemma \ref{Lem_Herman}.
\begin{lemma}[Epstein]
There cannot be any additional non-repelling periodic points.
\end{lemma}
\begin{proof}
Suppose that $h_\a$ has an additional non-repelling cycle
\[
\{z_1\mapsto z_2\mapsto \ldots, \mapsto z_k\mapsto z_1\}\;.
\]
Let $v_1$ and $v_2$ be the two critical values of $h_\a$, set
\[
X = \{0,1, z_1,\ldots,z_k\}, \quad X' = X\cup \{v_1,v_2\}\; .
\]
Let ${\cal Q}^1(X)$ (resp. ${\cal Q}^1(X')$) be the set of
meromorphic quadratic differentials on $\Cc$ which are holomorphic
outside $X$ (resp.\ $X'$) and have at most simple poles in $X$
(resp.\ $X'$). Let ${\cal Q}^2(X)$ be the set of meromorphic
quadratic differentials on $\Cc$ which are holomorphic outside $X$,
have at most double poles in $X$ and whose polar part of order $2$
along $X$ is of the form
\[
A \frac{dz^2}{z^2} + B \frac{dz^2}{(z-1)^2} + C\sum_{i=1}^k
\frac{dz^2}{(z-z_i)^2}\quad\text{with }A,B,C\in \C\;.
\]
The sets ${\cal Q}^1(X)$, ${\cal Q}^1(X')$ and ${\cal Q}^2(X)$ are
vector spaces of respective dimensions $k-3$, $k-1$ and $k$. We can
define a linear map $\nabla : {\cal Q}^2(X)\to {\cal Q}^1(X')$ as
follows. If $U$ is a simply connected subset of $\Cc\setminus X'$,
then $h_\a:h_\a^{-1}(U)\to U$ is a (trivial) covering map. We let
$(g_i:U\to  \Cc)_{i\in I}$ be the countably many inverse branches
and we set
\[
(h_\a)_* q |_U = \left(\sum_{i\in I} g_i^* q\right)\; .
\]
The sum is convergent because
\[
\sum_{i\in I} \int_U |g_i^* q| = \int_{h_\a^{-1}(U)} |q| <\infty\; .
\]
We can define in such a way a quadratic differential $(h_\a)_* q$
which is holomorphic outside $X'$. A local analysis shows that
\[
\nabla q := (h_\a)_* q-q
\]
has at most simple poles at points of $X'$ and thus, belongs to
${\cal Q}^1(X')$.

Since the dimension of ${\cal Q}^1(X')$ is less than the dimension
of ${\cal Q}^2(X)$, the linear map $\nabla$ is not injective and
there is a $q\in {\cal Q}^2(X)$ such that $\nabla q=0$, i.e.,
$(h_\a)_* q = q$. To see that this is not possible, set
\[
U_\e := D(0,\e)\cup D(1,\e) \cup \bigcup_{i=1}^k
h_\a^{-i}\bigl(D(z_1,\e)\bigr)\; ,\quad V_\e := h_\a^{-1}(U_\e)\;,
\]
let $W_\e\subset \Cc\setminus\bigl(U_\e\cup  \{v_1,v_2\}\bigr)$ be a
simply connected subset of full measure and let $g_i:W_\e\to \Cc$ be
the countably many inverse branches of $h_\a$. Then, for $\e$
sufficiently small, we have
\[
\int_{\Cc\setminus U_\e} \bigl|(h_\a)_*q\bigr| =  \int_{W_\e}
\left|\sum_{i}g_i^* q\right| \leq \sum_{i} \int_{W_\e} \bigl|g_i^*
q\bigr| = \int_{\C\setminus V_\e} |q|
\]
with equality if and only if each $g_i^* q$ is a (real positive)
multiple of $(h_\a)_*q=q$. In particular $q=h_\a^*(g_i^*q)$ has to
be locally, and thus globally, a constant multiple of $h_\a^*q$,
i.e.\ $q=c\cdot h_\a^* q$ for some constant $c> 0$. But in that case
$g_i^*q = c\cdot q$ and the sum $\displaystyle \sum_{i} \int_{W_\e}
\bigl|g_i^* q\bigr|$ will be diverging which is not the case. Thus,
\[
\int_{\Cc\setminus U_\e} |q| \leq  \int_{\Cc\setminus V_\e}
|q|-C_\e\quad\text{with }C_\e>0\; .
\]
Note that for $\delta<\e$, we have
\begin{eqnarray*}
\int_{\Cc\setminus U_{\delta}} |q| & = & \int_{\Cc\setminus U_{\e}}
|q| +\int_{U_{\e}\setminus
U_{\delta}} |q| \\
& \leq & \int_{\Cc\setminus V_\e} |q|-C_\e +
\int_{V_\e\setminus V_{\delta}} |q| \\
& = & \int_{\Cc\setminus V_{\delta}} |q| -C_\e\;,
\end{eqnarray*}
thus
\[\int_{\Cc\setminus V_{\delta}} |q| - \int_{\Cc\setminus
U_{\delta}} |q|\geq C_\e>0\;.
\]
We will obtain a contradiction by proving
\[\liminf_{\delta\to 0} \left(\int_{\Cc\setminus V_{\delta}} |q| -
\int_{\Cc\setminus U_{\delta}} |q| \right)\leq 0\;.
\]
This is the place where we use the fact that the cycle is
non-repelling. As $\delta\to 0$, we can find a radius $r_{\delta} =
\delta+o(\delta)$ such that
\[
D(0,r_{\delta})\cup D(1,r_{\delta}) \cup D(z_1,r_{\delta})\cup
\bigcup_{i=2}^k h_\a^{-i}\bigl(D(z_1,\delta)\bigr) \subset
V_{\delta}.\] Then, $U_{\delta}\setminus V_{\delta}$ is contained
within the union of three annuli
\[
\{z~;~r_{\delta}\leq |z|<\delta\}\cup \{z~;~r_{\delta}\leq
|z-1|<\delta\} \cup \{z~;~r_{\delta}\leq |z-z_1|<\delta\}\; .
\]
Since $q$ has at most double poles at $0$, $1$ and $z_1$, the
integral of $|q|$ on those annuli tends to $0$ as $\delta$ tends to
$0$ and we have
\[
\int_{\Cc\setminus V_{\delta}} |q| - \int_{\Cc\setminus U_{\delta}}
|q| = \int_{U_{\delta}\setminus V_{\delta}} |q| -
\int_{V_{\delta}\setminus U_{\delta}} |q| \leq
\int_{U_{\delta}\setminus V_{\delta}} |q|\underset{\delta\to
0}\longrightarrow 0\; .
\]
\end{proof}

\begin{lemma}
\label{Lem_Herman} There cannot be any cycle of Herman rings.
\end{lemma}
\begin{proof}
Recall that $0$ is a multiple fixed point and its immediate basin of
attraction must contain a critical point $\omega_0$ and the critical
value $v_0=h_\a(\omega_0)$. Also, $1$ is a Cremer point. It must be
accumulated by the orbit of the second critical point $\omega_1$
with critical value $v_1=h_\a(\omega_1)$.

Assume there is a cycle of Herman rings $H_1\mapsto H_2\mapsto
\ldots\mapsto H_k\mapsto H_1$. Let $\Gamma$ be the union of the
equators of the Herman rings $H_i$ ($\Gamma$ is the union of a cycle
of Jordan curves). Choose a connected component $W$ of $\Cc\setminus
\Gamma$ which does not contain $1$. Then, there are infinitely many
iterates of $v_1$ contained in $W$ (accumulating a boundary
component of some Herman ring). In particular, there is an integer
$m>2$ such that $h_\a^{\circ m}(v_1)\in W$. Let $D$ be a disk around
$1$ avoiding $\Gamma$, the forward orbit of $v_0$ and the $m$ first
iterates of $v_1$. Let $D_{-1}$ be the connected component of
$h_\a^{-1}(D)$ containing $1$. Since $D\setminus \{1\}$ does not
contain any singular value of $h_\a$, $h_\a:D_{-1}\to D$ has to be
an isomorphism. Since $D_{-1}$ contains $1$ and avoids $\Gamma$, it
does not contain $h_\a^{\circ m}(v_1)$. So, $D_{-1}$ is a disk
avoiding $\Gamma$, the forward orbit of $v_0$ and the $m$ first
iterates of $v_1$. We can therefore construct inductively a sequence
of disks $D_{-k}$ containing $1$ such that $h_\a^{\circ k}:D_{-k}\to
D$ is an isomorphism. Since $|(h_\a^{\circ k})'(1)|=1$ for all
$k\in\N$, by Koebe's one quarter theorem the disks $D_{-k}$ contain
a common neighborhood of $1$ on which the iterates of $h_\a$ form a
normal family. This contradicts the fact that $1$ is a Cremer point
contained in the Julia set of $h_\a$.
\end{proof}

Note that if we choose $\a\in \Q$, $N_\a$ will have a wandering
domain that projects to a parabolic basin of a parabolic fixed
point. If $\a$ is a Brjuno number, $N_\a$ will have a univalent
Baker domain of parabolic type II which projects to a Siegel disk of
$g_\a$.

\medskip
We can construct other examples in a similar way. The maps we will
present do not have fixed points. It follows from Proposition
\ref{Prop_NewtonMaps} that they are Newton maps of non-vanishing
entire functions, whose singularities over $0$ are therefore direct.

Assume
\[
        N(Z) = Z + \frac{\a}{1+\e\sin(2\pi Z)}
\]
with
\[0<\e<1 \quad\text{and}\quad 0<\a< m_\e =
\left\lfloor \frac{(1-\e)^2}{2\pi \e}\right\rfloor  \;.
\]
Then, $N$ is the Newton map of an entire function $f$ such that
$f(t)\to 0$ as $t\in \R\to +\infty$. The restriction of $N$ to $\R$
is an increasing homeomorphism which commutes with translation by
$1$. Indeed,
\[
N'(Z) = 1-\frac{2\pi \e \a \cos (2\pi Z)}{\bigl(1+\e \sin (2\pi
Z)\bigr)^2} \geq  1-\frac{2\pi \e \a}{(1-\e)^2} >0.\]
 Thus, it has a well defined rotation number ${\rm Rot}(N)$.
This rotation number is positive since $N(Z)>Z$. Note that for
$\a=m_\e$, $N(0) = m_\e$ and thus, ${\rm Rot}(N)= m_\e$. For each
fixed $\e\in (0,1)$, the rotation number increases continuously from
$0$ to $m_\e$ as $\a$ increases from $0$ to $m_\e$. If ${\rm
Rot}(N)$ is rational, then $N$ has a chain of wandering domains
along the real axis. If ${\rm Rot}(N)$ is a Brjuno number, $N$ has a
univalent Baker domain of hyperbolic type centered on the real axis.
For suitably chosen parameters $\a$, ${\rm Rot}(N)$ is irrational
and the induced map $N:\R/\Z\to \R/\Z$ is topologically but not
analytically conjugate to the rotation $Z\mapsto Z+{\rm
Rot}(N):\R/\Z\to \R/\Z$. It should follow that $N$ does not have any
Baker domain associated to the singularity of $f$ containing the
large positive real numbers. The proof should be similar to the one
we presented above: study the dynamics modulo $1$.

In the previous examples, $f$ had a direct singularity containing
critical points of $f$. One may wonder whether it is the presence of
critical points that prevents $N_f$ from having a Baker domain
associated to the singularity. The following example shows that this
is not the case. We still assume $\a>0$ and set
\[
        N_\a(Z)= Z+\a e^{e^{2i\pi Z}}\;.
\]
Then, $N_\a$ does not have any fixed points. So, it is the Newton
map of the non-vanishing entire function
\[
f_\a(Z) = \exp\left(-\frac{1}{\a}\int_0^Z e^{-e^{2i\pi W}}
dW\right)\;.
\]
Note that when $W\in \R$, the real part of $e^{-e^{2i\pi W}}$ is
greater than $1/e$. Thus, for $\a>0$ and for $t\in [0,+\infty)$, we
have
\[
|f_\a(t)| \leq e^{-t/(e\a)}\underset{t\to +\infty}\longrightarrow 0.
\]
The entire map $f_\a$ has a singularity over $0$ containing large
real numbers. This is a direct singularity since $f_\a$ does not
vanish. In addition, $N_\a$ does not have poles and so, $f_\a$ does
not have critical points.

Again, $N_\a(Z+1)=N_\a(Z)+1$ and $N_\a$ projects via $Z\mapsto
z=e^{2i\pi Z}$ to an entire map $g_\a$ fixing $0$ with multiplier
$e^{2i\pi \a}$:
\[
g_\a(z) = ze^{2i\pi\a e^z}\;.
\]
By a result of Bergweiler \cite{Bergweiler2}, the Fatou sets of
$N_\a$ and $g_\a$ correspond under the map $Z\mapsto e^{2i\pi Z}$.
If $g_\a$ has a Siegel disk around $0$, the map $N_\a$ has a Baker
domain of parabolic type II which corresponds to the singularity of
$f_\a$ described above. But if $g_\a$ has a Cremer point at $0$,
there can be no Baker domain for $N_\a$ associated to the
singularity of $f_\a$ described above.

\section{A Virtual Immediate Basin Implies an Asymptotic Value}
\label{Sec_Converse}

\begin{theorem}[Virtual Immediate Basin Contains Asymptotic Path]
\label{Thm_Converse} Let $f:\C\to\C$ be an entire function such that
its Newton map $N_f$ has a virtual immediate basin $V$. If $V$ is
parabolic of type I or type II, then $0$ is an asymptotic value of
$f$ with asymptotic path in $V$. There exists $H>0$ such that the
same is true if $V$ is hyperbolic with constant $h\geq H$.
\end{theorem}
Bergweiler, Drasin and Langley have constructed an
entire function for which $0$ is not an asymptotic value and whose
Newton map has a virtual immediate basin of hyperbolic type \cite{BDL}. Thus,
the statement of Theorem \ref{Thm_Converse} cannot be extended to
all hyperbolic virtual immediate basins.

Using Theorem \ref{Thm_Converse}, we can give the following
formulation of Theorem 5.1 in \cite{RS}.
\begin{corollary}[Outside Immediate Basins]
Let $N_f$ be the Newton map of an entire function $f$ and $U_\xi$
the immediate basin of the attracting fixed point $\xi\in \C$ for
$N_f$. Let $\Gamma_1,\Gamma_2\subset U_\xi$ be two $N_f$-invariant
curves connecting $\xi$ to $\infty$ such that $\Gamma_1$ and
$\Gamma_2$ are non-homotopic in $U_{\xi}$ and let $\tilde{V}$ be an
unbounded component of $\C\setminus (\Gamma_1\cup\Gamma_2)$. If the
set $N_f^{-1}(\{z\})\cap\tilde{V}$ is finite for all $z\in\Cc$, then
$f|_{\tilde{V}}$ assumes the value $0$ or has $0$ as an asymptotic
value.
\end{corollary}
\begin{proof} If $0\not\in f(\tilde{V})$, then the virtual immediate basin constructed in the proof of \cite[Theorem 5.1]{RS} is parabolic of type I.
\end{proof}
For the proof of Theorem \ref{Thm_Converse}, we will need the
following corollary to the Koebe distortion theorem. We thank Dierk
Schleicher for pointing it out to us.
\begin{lemma}[Bounded Non-Linearity]
\label{Cor_Koebe} Let $R>0$, $g:B_R(0)\to\C$ be univalent and
$\eps>0$. If $r/R$ is sufficiently small, then
\[
    \left| \frac{g(z)-g(w)}{g'(z)(z-w)}-1\right|< \eps
\]
for all $w,z\in B_r(0)$.
\end{lemma}
\begin{proof}
By possibly conjugating $g$ with $z\mapsto Rz$, multiplying $g$ with
a constant or adding a constant to $g$, we may assume that $R=1$,
$g(0)=0$ and $g'(0)=1$. Fix $0<r<1$. By the Koebe distortion
theorem, there is an $\alpha>0$ independent of $g$ such that
\[
|g(z)-g(w)-(z-w)g'(z)|<\alpha|(z-w)|^2
\]
for all $z,w\in B_r(0)$ (Taylor expansion around $z$). Moreover,
there is a $\beta>0$ so that $|g'(z)|>\beta$ for all $z\in B_r(0)$.
This yields
\[
\left|\frac{g(z)-g(w)}{g'(z)(z-w)}-1\right|<\alpha\left|\frac{z-w}{g'(z)}\right|
<\frac{2\alpha r}{\beta} \;.
\]
It follows from the Koebe distortion theorem that $\alpha\to 0$ and
$\beta\to 1$ as $r\to 0$. The claim follows.
\end{proof}
\proofof{Theorem \ref{Thm_Converse}} Suppose first that $V$ is
parabolic of type I. Then, there exists a weakly absorbing set $A$ of $V$ and a conformal conjugacy
$(\C,\phi,T)$ such that $F:=\phi(A)$ is an absorbing set for $T:z\mapsto z+1$ in $\C$. Since $\phi|_A$ is univalent, it has a univalent inverse $\psi:F\to A$. With this, we get for $z\in F$
that $N_f(\psi(z)) = \psi(z+1)$, and hence
\[
    \psi(z)-\frac{f(\psi(z))}{f'(\psi(z))}=\psi(z+1)\;.
\]
It follows that
\begin{equation}
\label{Eqn_Estimate0}
    \frac{f'(\psi(z))}{f(\psi(z))}\cdot \big(\psi(z+1)-\psi(z)\big) = -1
\end{equation}
(note that since $V$ is a virtual immediate basin, $f$ has no roots
in $\psi(F)$). Let $0 < \eps <1$. By Lemma \ref{Cor_Koebe}, there exists
$R>2$ such that if $B_R(z) \subset F$, then
\begin{equation}
\label{Eqn_Estimate1}
    \left|\frac{\psi'(z)}{\psi(z+1)-\psi(z)}-1\right|<\eps\;,
\end{equation}
and by equation (\ref{Eqn_Estimate0}) and inequality
(\ref{Eqn_Estimate1}) we get
\begin{equation}
\label{Eqn_Estimate2}
\left|\frac{f'(\psi(z))}{f(\psi(z))}\cdot\psi'(z)+1
\right|=\left|\frac{f'(\psi(z))}{f(\psi(z))}\cdot\psi'(z)\cdot\frac{\psi(z+1)-\psi(z)}{\psi(z+1)-\psi(z)}+1
\right|<\eps\;.
\end{equation}
Since $F$ contains all sufficiently far right translates of the disk $B_R(z_0)$, for every $z_0\in F$ there exists $S_{z_0}\geq 0$ such that 
(\ref{Eqn_Estimate2}) holds for all $z_0+t$ with real $t\geq S_{z_0}$. 

Let $z_0\in F$ such that $S_{z_0}=0$. Then, for $t\geq 0$ and $z=z_0+t\in F$, we use a standard estimate in complex variables and inequality
(\ref{Eqn_Estimate2}) to get
\begin{eqnarray*}
    \left| \log(f(\psi(z)))+z \right| &\leq& \left|\int_{z_0}^z \left((\log\circ f\circ\psi)'(\zeta)+1\right) d\zeta \right| + \left| \log(f(\psi(z_0)))+z_0\right| \\
&\leq&\sup_{w\in [z_0,z]} \left\{  \left| \frac{f'(\psi(w))}{f(\psi(w))}\cdot \psi'(w)+1\right|\right\} \cdot |z-z_0| + C'\\
&\leq&\eps\cdot|z-z_0|+C' \\
&\leq& \eps\cdot |z|+C\;,
\end{eqnarray*}
where $C'=|\log(f(\psi(z_0)))+z_0|$ and $C>0$ depend only on $z_0$;
$[z_0,z]$ denotes the straight line segment in $F$ connecting
$z_0$ to $z$. It follows that $\log(f(\psi(z)))\in B_{\eps
|z|+C}(-z)$ and 
\begin{equation}
\label{Eqn_Estimate3}
    \Re(\log(f(\psi(z)))) < - \Re(z) + \eps |z|+C\;.
\end{equation}
Since $\Im(z)$ does not depend on $t$, we have that $|z|/\Re(z)\to 1$ as $t\to\infty$ and the right hand side of inequality (\ref{Eqn_Estimate3}) converges to $-\infty$. Hence, exponentiating (\ref{Eqn_Estimate3}) yields
$\lim_{t\to+\infty}f(\psi(z))=0$.

Analogous estimates hold for sufficiently large imaginary parts if
$V$ is parabolic of type II. If $V$ is hyperbolic, sufficiently large $h$ will permit a construction as above. This finishes the proof. 
\qed
\remark In fact, we not only show the existence of an asymptotic
path to $0$ for $f$ in $V$, but even that $V$ has an $N_f$-invariant
open subset in which $f$ converges to $0$ along $N_f$-orbits. This
is another similarity between immediate basins and their virtual
counterparts.

\section{Acknowledgements}
We thank Adrien Douady for raising the question of a relation
between virtual immediate basins and asymptotic values and Dierk
Schleicher for his helpful comments and his support. We also thank
Walter Bergweiler and Alexandre Eremenko for several interesting
discussions in which we learned a lot about transcendental
functions.

\end{document}